\documentclass[12pt]{amsart}

\usepackage{amsmath,amsthm,amssymb}
\font\teneufm=eufm10
\font\seveneufm=eufm7
\font\fiveeufm=eufm5
\newfam\frakturfam

\textfont\frakturfam=\teneufm
\scriptfont\frakturfam=\seveneufm
\scriptscriptfont\frakturfam=\fiveeufm

\newtheorem{lm}{Lemma}
\newtheorem{theor}{Theorem}

\newtheorem{rem}{Remark}

\newtheorem{prob}{Problem}
\def\bee{\begin{eqnarray}}
\def\bes{\begin{eqnarray*}}
\def\eee{\end{eqnarray}}
\def\ees{\end{eqnarray*}}

\def\a{\alpha}

\def\s{\sigma}

\def\Proof{{\sl Proof.}\ }

\title{Erroneous proofs of the wildness of some automorphisms of free metabelian Lie algebras}

\begin{document}
\date{}
\maketitle

\begin{center}
{\bf Ualbai Umirbaev}\footnote{Department of Mathematics,
 Wayne State University,
Detroit, MI 48202, USA 
and Institute of Mathematics and Mathematical Modeling, Almaty, 050010, Kazakhstan
e-mail: {\em umirbaev@wayne.edu}}
\end{center}

\begin{abstract} 
The well known Bachmuth-Mochizuki-Roman'kov Theorem \cite{BM,Romankov85} states that every automorphism of the free metabelian group of rank $\geq 4$ is tame. In 1992 Yu. Bahturin and S. Nabiyev \cite{BN} claimed that every nontrivial inner automorphism of the free metabelian Lie algebra $M_n$ of any rank $n\geq 2$ over a field of characteristic zero is wild.  More examples of wild automorphisms of $M_n$ of rank $n\geq 4$ were given in 2008 by Z. \"Ozcurt and N. Ekici \cite{OE}.
The main goal of this note is to show that both articles contain uncorrectable errors and to draw the attention of specialists to the fact that the question of tame and wild automorphisms for free metabelian Lie algebras $M_n$ of rank $n\geq 4$ is still widely open. 
\end{abstract}



\section{Introduction}

\hspace*{\parindent}

Let $K$ be an arbitrary field, and let ${\mathfrak M}$ be an arbitrary homogeneous variety of linear algebras over $K$. Denote by $A=K_{\mathfrak M}\langle x_1,x_2,\ldots,x_n\rangle$
  the free algebra of ${\mathfrak M}$ with a free set of generators $X=\{x_1,x_2,\ldots,x_n\}$. If ${\mathfrak M}$ is a unitary variety of algebras (see, for example \cite{KBKA}) then we assume that $A$ has an identity element $1$.  
 Let $\mathrm{Aut}(A)$ be the group of all automorphisms of $A$. Denote by $\phi = (f_1,f_2,\ldots,f_n)$ the automorphism 
  $\phi$ of $A$ such that $\phi(x_i)=f_i,\, 1\leq i\leq n$. An automorphism  
\bes
(x_1,\ldots,x_{i-1}, \a x_i+f, x_{i+1},\ldots,x_n),
\ees
where $0\neq\a\in K,\ f\in K_{\mathfrak M}\langle X\setminus \{x_i\}\rangle$, 
 is called {\em elementary}.    The subgroup $\mathrm{TAut}(A)$ of $\mathrm{Aut}(A)$ generated by all 
  elementary automorphisms is called the {\em tame automorphism group}, 
 and the elements of this subgroup are called {\em tame automorphisms} 
 of $A$. Nontame automorphisms of $A$ are called {\em wild}.

The question about tame and wild automorphisms of free algebras plays an important role in combinatorial algebra. 
It is well known that all automorphisms of the polynomial algebra $K[x,y]$ in two variables $x,y$ over a field $K$ are tame \cite{Jung, Kulk}. Similar results hold for free associative algebras \cite{Czer, ML70} and for free Poisson algebras (in characteristic zero) \cite{MLTU} (see also \cite{MLSh,MLU16}).  Moreover, the automorphism groups of polynomial algebras, free associative algebras, and free Poisson algebras in two variables are isomorphic. 

The automorphism  groups of commutative and associative algebras generated by three elements are much more complicated. 
The well-known Nagata automorphism \cite{Nagata} 
\bes
\sigma=(x+2y(zx-y^2)+z(zx-y^2)^2, y+z(zx-y^2),z)
\ees
of the polynomial algebra $K[x,y,z]$ over a field $K$ of characteristic $0$ is proven to be wild \cite{Umi25}. The so-called Anick automorphism (see \cite[p. 343]{Cohn})
\bes
\delta&=&(x+z(xz-zy),\,y+(xz-zy)z,\,z)
\ees
of the free associative algebra $K\langle x,y,z\rangle$ over a field $K$ of characteristic $0$ is also proven to be wild \cite{Umi33}. 
The Nagata
automorphism gives an example of a wild
automorphism of free Poisson algebras in three variables. In 2006 I.P. Shestakov (unpublished) constructed an analogue of the Anick automorphism 
\bes
\tau=(x+\{xz-\{y,z\},z\}, y + z(xz-\{y,z\}), z)
\ees  
for free Poisson algebras in three variables. 

 It is well known \cite{Smith} that the Nagata and the Anick automorphisms are stably tame.  

In 1981 Shafarevich \cite{Shafarevich81} and in 1983 Anick \cite{Anick} independently proved that every automorphism of the polynomial algebra $K[x_1,x_2,\ldots,x_n]$ over a field $K$ of characteristic zero can be approximated by tame automorphisms with respect two different but similar topologies. This can be considered as a weaker generalization of the Jung-van der Kulk Theorem for polynomial algebras of rank $n\geq 3$.  In particular, the Nagata automorphism is wild but can be approximated by tame automorphisms. 

 In 1964 P. Cohn  proved \cite{Cohn} that all automorphisms of finitely generated free
Lie algebras over a field are tame.
 Later this result was extended to free algebras of Nielsen-Schreier varieties \cite{Lewin}. Recall that
a variety of universal algebras is called Nielsen-Schreier, if any subalgebra of a free algebra of this variety is free,
i.e., an analog of the classical Nielsen-Schreier theorem is true.
 The varieties of all non-associative algebras \cite{Kurosh},
 commutative and anti-commutative algebras \cite{Shirshov54}, Lie algebras and Lie $p$-algebras 
 \cite{Shirshov53,Witt}, and Lie superalgebras  \cite{Mikhalev85,Stern} and Lie $p$-superalgebras \cite{Mikhalev88} over a field are Nielsen-Schreier. 
 Some other examples of Nielsen-Schreier 
 varieties can be found in \cite{SU02,MS14,U94,U96,Chibrikov}. 
 
It was recently shown \cite{DU22} that the varieties of pre-Lie (also known as right-symmetric) algebras and Lie-admissible algebras over a field of characteristic zero are Nielsen--Schreier. In particular, every automorphism of a free right-symmetric and a free Lie-admissable algebra of finite rank is tame. The tameness of automorphisms of free right-symmetric algebras in two variables over any field was proven in \cite{KMLU}

One of the classical and well studied varieties of algebras is the variety of metabelian Lie algebras. Let $L_n$ be a free Lie algebra of rank $n$ in the variables $x_1,\ldots,x_n$ over a field $K$. 
Then $M_n = L_n/L''_n=L_n/[[L_n,L_n],[L_n,L_n]]$
 is the free metabelian Lie algebra 
of rank $n$ in the variables $y_i=x_i+L''_n$, where $1\leq i\leq n$. 
Let $0\neq z\in [M_n,M_n]$. Then $\mathrm{ad}\,z : M_n\to M_n (m\mapsto [z,m])$ is a locally nilpotent derivation and $(\mathrm{ad}\,z)^2=0$. The {\em exponential} 
$\exp(\mathrm{ad}\,z)$ of $M_n$ can be written as  
\bes
\exp(\mathrm{ad}\,z)=\mathrm{id}+\mathrm{ad}\,z =(y_1+[z,y_1],\ldots,y_n+[z,y_n]). 
\ees
Automorphisms of this type are called {\em inner}. 

Obviously, the group of tame automorphisms of $M_2$ coincides with the group of linear automorphisms $\mathrm{GL}_2(K)$. It follows that any nontrivial inner automorphism of $M_2$ is wild (see also \cite[Proposition 4]{Shmelkin} and \cite[Theorem 3]{Papistas}). 
In 1992, V. Drensky \cite{Drensky92} proved that the inner automorphism $\mathrm{exp}(\mathrm{ad}[x_1,x_2])$ of $M_3$ is wild. In 2008, V. Roman'kov \cite{Romankov08} proved that the group of automorphisms of $M_3$ cannot be finitely generated modulo the subgroup generated by all inner and tame automorphisms. 

An analogue of the Anick-Shafarevich Theorem \cite{Anick,Shafarevich81} for free metabelian Lie algebras was proven in 1993 by Bryant and Drensky \cite{BD93-2}. They proved that the group of tame automorphisms is dense in the group of all automorphisms of $M_n$ if $n\geq 4$. A sharpened version of this result can be found in \cite{KP16,Nauryzbaev}. 

In 1992 Yu. Bahturin and S. Nabiyev \cite{BN} announced that any nontrivial exponential automorphism of the free metabelian Lie algebra $M_n$ of rank $n\geq 2$ over a field of characteristic zero is wild. Unfortunately the proof contains some uncorrectable mistakes. The paper \cite{BN} is not well written and it seems that many specialists in the area accepted the main result without checking the details of the proof. We explain the main steps of the proof given in \cite{BN} using more modern terminology and explicitly show the first main mistake in Section 2.

One more example of a wild automorphism of $M_n$ for $n\geq 4$ was given by Z. \"Ozcurt and N. Ekici \cite{OE}. The proof given in \cite{OE}, in fact, contradicts some results from \cite{BD93-1,BD93-2}. We discuss this in Section 3. 

  Thus the following problem is wide open. 
\begin{prob}\label{pr1} Is any automorphism of the free metabelian Lie algebra $M_n$ of rank $n\geq 4$ over any field tame? 
\end{prob}

The well known Bachmuth-Mochizuki-Roman'kov Theorem \cite{BM,Romankov85} states that every automorphism of the free metabelian group of rank $\geq 4$ is tame. More information about automorphisms of free metabelian groups can be found in \cite{Romankov20}.

In 1992 V. Drensky \cite{Drensky92}, in 1993 Papistas 
\cite{Papistas93}, and in 1995 Bahturin and Shpilrain \cite{BSh95} proved that free algebras of rank $n\geq 2$ of any polynilpotent variety of Lie algebras that is not abelian or metabelian have wild automorphisms. This result also makes the case of free metabelian algebras more important.

\section{On the paper by Bahturin and Nabiyev}

\hspace*{\parindent}

\subsection{The structure of $M_n$} 

\hspace*{\parindent}

Let $M_n$ be the free metabelian Lie algebra over a field $K$ freely generated generated by $x_1,\ldots,x_n$. We give a description of $M_n$ using Shmel'kin's wreath products \cite{Shmelkin}. Let $Y_n$ be the abelian Lie algebra with a linear basis $y_1,y_2,\ldots,y_n$. The universal enveloping algebra $U(Y_n)$ is the polynomial algebra $U(Y_n)=K[y_1,\ldots,y_n]$. Let $T_n$ be the free (left) $U(Y_n)$-module with a basis $t_1,\ldots,t_n$. Consider the direct sum 
\bes
M=Y_n\oplus T_n. 
\ees
 Turn $M$ into a Lie algebra by 
\bes
[a+t,b+s]= as-bt 
\ees
where $a,b\in Y_n$ and $t, s\in T_n$. Then the subalgebra of $M$ generated by 
\bes
x_i=y_i+t_i, \ \ 1\leq i\leq n, 
\ees
is the free metabelian Lie algebra $M_n$ with the free variables $x_1,\ldots,x_n$ \cite{Shmelkin}. 

\subsection{The Jacobian matrix of an endomorphism of $M_n$}

\hspace*{\parindent}

Let $f\in M_n$ be an arbitrary element $M_n$. Then $f$ is uniquely represented as 
\bes
f=y+t, y\in Y_n, t\in T_n. 
\ees
Moreover, $t$ is uniquely represented as 
\bes
t=d_1t_1+\ldots+d_nt_n, d_1,\ldots,d_n\in K[y_1,\ldots,y_n]=U(Y_n). 
\ees
We call $d_i=\frac{\partial f}{\partial x_i}$ the Fox derivatives of $f$  (see \cite{U93}). Set  
\bes
\partial(f)=(d_1,\ldots,d_n)=\big(\frac{\partial f}{\partial x_1},\ldots,\frac{\partial f}{\partial x_n}\big). 
\ees
Also set $Y=(y_1,\ldots,y_n)^t$, where $t$ means transpose. There is a well known statement (\cite[Lemma 5]{U93}) that says $f\in [M_n,M_n]$ if and only if 
\bee\label{f1}
\partial(f)Y=0. 
\eee

Let $\phi=(f_1,\ldots,f_n)$ be an arbitrary endomorphism of $M_n$ such that $\phi(x_i)=f_i$ for all $i$. In particilar, $X=(x_1,\ldots,x_n)$ is the identity automorphism. We define the Jacobian matrix of $\phi$ by 
\bes
J(\phi)=[\frac{\partial f_i}{\partial x_j}]_{1\leq i,j\leq n}=\left[\begin{array}{ccc}
 \partial(f_1)\\
 \ldots\\
 \partial(f_n)\\
\end{array}\right].
\ees
If $\phi=(f_1,\ldots,f_n)$ and $\psi=(g_1,\ldots,g_n)$ then 
\bes
\phi\circ\psi(x_i)=g_i(f_1,\ldots,f_n)
\ees
for all $1\leq i\leq n$. The corresponding Chain Rule gives \cite{U93} that 
\bee\label{f2}
J(\phi\circ\psi)=\overline{\phi}(J(\psi)) J(\phi), 
\eee
where $\overline{\phi}$ denotes the endomorphism of $K[y_1,\ldots,y_n]$, first induced by $\phi$ to $Y_n$ and then extended to $U(Y_n)$. 

Denote by $\mathrm{IAut}\,M_n$ the group of all automorphisms of $M_n$ that induces the identity automorphism on $M_n/[M_n,M_n]$. 
\begin{lm}\label{l1}
 The map 
\bes
J : \mathrm{IAut}(M_n) \to \mathrm{GL}_n(U(Y_n))
\ees
is an injective antimorphism of groups. 
\end{lm}

This is a direct corollary of (\ref{f2}) since $\overline{\phi}=\mathrm{id}$. 

\begin{rem} The matrix $E+QA_n$ in \cite[Theorem 2]{BN} is the Jacobian matrix of an automorphism. 
It is mentioned (without proof) in the middle of \cite[Theorem 2]{BN} that $J$ is an isomorphism. Sometimes it depends on the definition of composition of automorphisms.  
\end{rem}

\subsection{Tame automorphisms}

\hspace*{\parindent}

The group of tame automorphisms $\mathrm{TAut}\,M_n$ is generated by all linear automorphisms and  the automorphisms of the form 
\bes
\phi_f= (x_1+f(x_2,\ldots,x_n), x_2,\ldots, x_n), f\in [M_n,M_n]. 
\ees 
Notice that 
\bes
J(\phi_f)=E+\left[\begin{array}{cccc}
 \partial(f)\\
0\\
 \ldots\\
 0\\
\end{array}\right]=E +e_1 \partial(f), 
\ees
where $E$ is the identity matrix and $e_1=(1,0,\ldots,0)^t$, since $\partial(f)=(0,\frac{\partial f}{\partial x_2},\ldots,\frac{\partial f}{\partial x_n})$. Notice that $\partial(f) e_1=0$. 

\begin{lm}\label{l2} Let $\a$ be an arbitrary linear automorphism of $M_n$. Then 
\bes
J(\a\phi_f\a^{-1})=E+\Phi\Psi, 
\ees
where $\Phi$ is a column vector and $\Psi$ is a row vector such that 
\bes
\Psi\Phi=0, \Psi Y=0. 
\ees
\end{lm}
\Proof Let $J(\a)=A$. Using (\ref{f2}) we get 
\bes
J(\a\phi_f\a^{-1})=\overline{\a}(J(\phi_f\a^{-1}) J(\a)=\overline{\a}(\overline{\phi_f}(J(\a^{-1}))J(\phi_f)) A\\
=\overline{\a}(A^{-1}J(\phi_f)) A=A^{-1}\overline{\a}(J(\phi_f)) A=A^{-1}\overline{\a}(E +e_1 \partial(f)) A\\
=A^{-1}(E +e_1 \overline{\a}(\partial(f))) A=E +A^{-1}e_1 \overline{\a}(\partial(f)) A=E+\Phi\Psi, 
\ees
where $\Phi=A^{-1}e_1$ and $\Psi= \overline{\a}(\partial(f)) A$. We have 
\bes
\Psi\Phi= \overline{\a}(\partial(f))e_1=\overline{\a}(\partial(f)e_1)=0. 
\ees
Obviously, 
\bes
\a\phi_f\a^{-1}=X+\a(f,0,\ldots,0)^t\a^{-1}=X+(f_1,\ldots,f_n)^t
\ees
for some $f_1,\ldots,f_n\in [M_n,M_n]$. Consequently, 
\bes
J(\a\phi_f\a^{-1})=E+\left[\begin{array}{cccc}
 \partial(f_1)\\
 \ldots\\
  \partial(f_n)\\
\end{array}\right]
\ees
and 
\bes
\Phi\Psi=\left[\begin{array}{ccc}
 \partial(f_1)\\
 \ldots\\
  \partial(f_n)\\
\end{array}\right]. 
\ees
By (\ref{f1}) we get 
\bes
\Phi\Psi Y=\left[\begin{array}{ccc}
 \partial(f_1)\\
 \ldots\\
  \partial(f_n)\\
\end{array}\right]Y=\left[\begin{array}{cccc}
 \partial(f_1)Y\\
 \ldots\\
  \partial(f_n)Y\\
\end{array}\right]=0. 
\ees

Notice that $\Psi Y$, the product of a row and a column, belongs to $U(Y_n)$. Since $\Phi=A^{-1}e_1$ is a nonzero column it follows that $\Psi Y=0$. $\Box$

\subsection{Inner automorphisms}

\hspace*{\parindent}

Let $0\neq z\in [M_n,M_n]$. Then $\mathrm{ad}\,z : M_n\to M_n (m\mapsto [z,m])$ is a locally nilpotent derivation and $(\mathrm{ad}\,z)^2=0$. The {\em exponential} automorphism of $M_n$  
\bes
\exp(\mathrm{ad}\,z)=\mathrm{id}+\mathrm{ad}\,z =(x_1+[z,x_1],\ldots,x_n+[z,x_n])^t. 
\ees
is called {\em inner}. 

Notice that $\partial([z,x_1])=-y_1\partial(z)$ since $z\in [M_n,M_n]\subseteq T_n$ and $[z,x_1]=[z,y_1+t_1]=-y_1 z$. Consequently, 
\bee\label{f3}
J(\exp(\mathrm{ad}\,z))=E-\left[\begin{array}{ccc}
y_1 \partial(z)\\
 \ldots\\
  y_n\partial(z)\\ 
\end{array}\right]= E-Y\partial(z).
\eee

\subsection{Wild automorphisms}

\hspace*{\parindent}

\begin{theor}\label{t1} \cite[Theorem 4]{BN}
 Every nontrivial inner automorphism $\s=\exp(\mathrm{ad}\,z)$ is wild. 
\end{theor}
\Proof Suppose that $\s$ is tame. Then $\s$ is a product of linear conjugates of the automorphisms of the type $\phi_f$ since $\s$ induces the identity automorphism on $Y_n$. Taking the Jacobian matrices, by Lemmas \ref{l1} and \ref{l2}, we get an equality of the form 
\bes
(E+\Phi_1\Psi_1)\ldots (E+\Phi_n\Psi_n)= E-Y\partial(z). 
\ees
From now on we deal only with the matrix multiplication. Recall that $\Phi_i$ are columns and $\Psi$ are rows. Let's consider the case when $n=3$ for clarity (the case $n=2$ is trivial). Then 
\bee\label{f4}
(E+\Phi_1\Psi_1)(E+\Phi_2\Psi_2)(E+\Phi_3\Psi_3)= E-Y\partial(z).
\eee
Set $\lambda_{ij}=\Psi_i\Phi_j\in K[y_1,\ldots,y_n]=U(Y_n)$.  Then (\ref{f4}) implies that 
\bee\label{f5}
\nonumber \Phi_1\Psi_1+\Phi_2\Psi_2+\Phi_3\Psi_3+\Phi_1\lambda_{12}\Psi_2+\Phi_1\lambda_{13}\Psi_3+\Phi_2\lambda_{23}\Psi_3\\
+\Phi_1\lambda_{12}\lambda_{23}\Psi_3=-Y\partial(z).
\eee
Notice that $\lambda_{ii}=0$ by Lemma \ref{l2}. Multiplying this equality by $\Psi_1$ from the left hand side, we get 
\bee\label{f6}
\lambda_{12}\Psi_2+\lambda_{13}\Psi_3+\lambda_{12}\lambda_{23}\Psi_3=0. 
\eee
This exactly corresponds to \cite[equation (17)]{BN}. 

\begin{rem} Notice that $\Psi_k=\overline{\Psi_k}A_n$ in \cite{BN} and the multiplication on the left first by $A_N$ and then by $\overline{\Psi_k}$ is equivalent to the multiplication by $\Psi_k$. 
\end{rem}

\subsection{The mistake}

\hspace*{\parindent}

Representation of (\ref{f4}) in the form 
\bes
(E+\Phi_1\Psi_1)(E+\Phi_2\Psi_2)(E-\Phi_1\Psi_1)(E+\Phi_1\Psi_1)(E+\Phi_3\Psi_3)\\
=(E+(\Phi_2+\lambda_{12}\Phi_1)(\Psi_2-\lambda_{21}\Psi_1))(E+\Phi_1\Psi_1)(E+\Phi_3\Psi_3)
= E-Y\partial(z)
\ees
is useless since the equation (\ref{f5}) corresponding to this equation is the same. 

But, interestingly, multiplication of (\ref{f5})
by $\Psi_2-\lambda_{21}\Psi_1$ leads to a "reduced" result. Of course, we can get the same result systematically, first multiplying  (\ref{f5}) by $\Psi_2$ and then subtracting $\lambda_{21}$ times  (\ref{f6}). In fact, multiplying  (\ref{f5}) by $\Psi_2$ gives
\bee\label{f7}
\lambda_{21}\Psi_1+\lambda_{23}\Psi_3+\lambda_{21}\lambda_{12}\Psi_2+\lambda_{21}\lambda_{13}\Psi_3
+\lambda_{21}\lambda_{12}\lambda_{23}\Psi_3=0.
\eee
Combining (\ref{f6}) and (\ref{f7}) we get 
\bee\label{f8}
\lambda_{21}\Psi_1+\lambda_{23}\Psi_3=0.
\eee
 We can also get this equality directly from (\ref{f5}) by multiplying it by $\Psi_2-\lambda_{21}\Psi_1$ and this equality corresponds to the last equality on \cite[page 51]{BN}. Unfortunately this equation contains the additional summand $\lambda_{21}\Psi_1$. 

This is the main mistake. It does not allow us to make further "triangular" calculations.

\section{On the paper by \"Ozkurt and Ekici}

\hspace*{\parindent}

The authors of \cite{OE} claim, in particular, that the automorphism 
\bes
\psi=(x_1+[[x_1,[x_2,x_3]],x_4], x_2,\ldots,x_n)
\ees
of $M_n$ is wild for all $n\geq 4$.

Let 
\bes
M_n=A_1\oplus \ldots \oplus A_k\oplus\ldots 
\ees
be the standard degree grading of $M_n$, i.e., 
$A_k$ is the linear span of (Lie) monomials of degree $k$. 

\bigskip 

Let $\mathrm{IAut}_i(M_n)$ be the subgroup of all automorphisms of $A$ that induces the identity automorphism on the factor-algebra 
\bes
M_n/(A_{i+1}+A_{i+2}+\ldots). 
\ees
We have 
\bes
\mathrm{Aut}(M_n)=\mathrm{IAut}_0(M_n)\supset \mathrm{IAut}_1(M_n)\supset \ldots \supset \mathrm{IAut}_k(M_n)\supset \ldots 
\ees 
Obviously, $\mathrm{IAut}(M_n)=\mathrm{IAut}_1(M_n)$. 

The automorphism $\psi$ belongs to $\mathrm{IAut}_3(M_n)$. By the main result of \cite{BD93-2} there exists a tame automorphism $\overline{\phi}$ of $M_n$ such that $\overline{\phi}\in \mathrm{IAut}_3(M_n)$ and $\overline{\phi}^{-1}\psi\in \mathrm{IAut}_4(M_n)$. This means that 
\bes
\overline{\phi}=(x_1+[[x_1,[x_2,x_3]],x_4]+u_1, x_2+u_2,\ldots,x_n+u_n), 
\ees
where $u_1,\ldots,u_n\in A_5+A_6+\ldots$. 

Let $L_n$ be the free Lie algebra of rank $n$ in the variables $z_1,\ldots,z_n$ such that $M_n = L_n/L''_n=L_n/[[L_n,L_n],[L_n,L_n]]$
 is the free metabelian Lie algebra 
of rank $n$ in the variables $x_i=z_i+L''_n$, where $1\leq i\leq n$. 

Suppose that $\overline{\phi}$ is induced by the tame automorphism $\phi$ of $L_n$. Then 
\bes
\phi=(z_1+[[z_1,[z_2,z_3]],z_4]+v_1+w_1, z_2+v_2+w_2,\ldots,z_n+v_n+w_n), 
\ees
where $v_1,\ldots,v_n\in L''_n$ are homogeneous polynomials of degree $4$ and $w_1,\ldots,w_n$ does not contain Lie monomials of degree less than $5$. 

As in \cite{OE}, denote by $\frac{\partial f}{\partial z_i}=\frac{\partial}{\partial z_i}(f)\in U(L_n)$ the (left) Fox derivatives of $f$ \cite{Reutenauer,Shpilrain,U93}, i.e., 
\bes
f=\frac{\partial f}{\partial z_1}z_1+\ldots+\frac{\partial f}{\partial z_n}z_n. 
\ees
Then 
\bes
\frac{\partial}{\partial z_1}([[z_1,[z_2,z_3]],z_4]+v_1)+\frac{\partial}{\partial z_2}(v_2)+\ldots+\frac{\partial}{\partial z_n}(v_n)\\
=z_4[z_2,z_3]+\frac{\partial}{\partial z_1}(v_1)+\frac{\partial}{\partial z_2}(v_2)+\ldots+\frac{\partial}{\partial z_n}(v_n)   \in [U(L_n),U(L_n)]
\ees
by \cite[Theorem 3.7]{BD93-1}. 

This contradicts the final conclusion of the proof of \cite[Theorem 4.1]{OE}.

\end{document}